# A FINITE-VOLUME DISCRETIZATION FOR DEFORMATION OF FRACTURED MEDIA

Eren Ucar[a], Eirik Keilegavlen[a], Inga Berre[a,b], Jan Martin Nordbotten[a,c]

[a] Department of Mathematics, University of Bergen, Bergen, Norway
[b] Christian Michelsen Research, Bergen, Norway
[c] Department of Civil and Environmental Engineering, Princeton University, Princeton, NJ, USA
Corresponding author: Eren Ucar, email: eren.ucar@uib.no

## ABSTRACT

Simulating the deformation of fractured media requires the coupling of different models for the deformation of fractures and the formation surrounding them. We consider a cell centered finite-volume approach, termed the multipoint stress approximation (MPSA) method, which is developed in order to discretize coupled flow and mechanical deformation in the subsurface. Within the MPSA framework, we consider fractures as co-dimension one inclusions in the domain, with the fracture surfaces represented as line pairs in 2D (faces in 3D) that displace relative to each other. Fracture deformation is coupled to that of the surrounding domain through internal boundary conditions. This approach is natural within the finite-volume framework, where tractions are defined on surfaces of the grid. The MPSA method is capable of modeling deformation considering open and closed fractures with complex and nonlinear relationships governing the displacements and tractions at the fracture surfaces. We validate our proposed approach using both problems for which analytical solutions are available and more complex benchmark problems, including comparison with a finite-element discretization.

## 1. INTRODUCTION

Subsurface rock is a porous medium containing fluids under complex in situ stress conditions [1,2]. Advances in the understanding of fluid flow and rock mechanics have vital importance in the development of several subsurface applications, including exploitation of geothermal-energy systems [3], enhanced recovery from oil and gas reservoirs [4], $CO_2$ storage [4,5] and underground storage of natural gas [6]. Deformation of the subsurface due to engineering operations or natural processes involves handling of structures with slit-like discontinuities (such as fractures) because of the constitution of subsurface rock. In many situations, these structures will dominate the mechanical behavior; thus, feasible modeling of rock mechanics for subsurface applications requires effective incorporation of fractures. The objective of this work is to describe and implement a numerical method to model mechanical deformations in fractured formations. Such a model is expected to be an important contribution to efficient hydromechanical coupling, particularly for subsurface applications.

In the literature on computational analysis of fractured domains, two main approaches have been applied: the boundary-element method (BEM) and the finite-element method (FEM) [7,8]. The BEM has been effectively applied to fracture problems for the past several decades, and its area of application has recently widened [3,9-13]. The BEM can accommodate problems including nonhomogeneous materials using Green's functions [9] and can be extended to 3D simulations [14]. With the FEM, the method must be 'enhanced' to accommodate discontinuous deformations [15]. In FEM simulations, classical approaches to incorporating fractures into the model are Lagrange multipliers [16,17], the penalty method [18], and the augmented Lagrangian method.

Strategies for multi-physics coupling, such as flow and mechanical deformation, have been analyzed and applied for several problems [10,19-21]. In most studies, whereas the mechanical equations are approximated by the FEM or BEM, the flow equations are generally discretized by conservative schemes such as the finite-volume method (FVM). The coupling of different numerical methods may require a different data structure for the disparate schemes. Moreover, it is also common to use two different software packages to solve the coupled problem, which requires an additional iteration procedure.

An alternative that avoids these difficulties is to discretize the full poro-mechanical system with the FVM. Specifically, we consider a cell-centered FVM for elastic deformation in porous media as a counterpart to finite-volume flow calculations in porous media [22]. The method is called the multi-point stress approximation (MPSA) due to its similarities with the multi-point flux approximation methods (e.g. [23]), which have been developed in the context of flow problems. The MPSA methods, and in particular the so-called MPSA-W method extended herein, can handle most polygonal and polyhedral grids, including both simplex and Cartesian-type grids, and have been shown to be robust with respect to material discontinuities [24,22]. The convergence properties of MPSA are established for elasticity [25], and poro-elasticity [26], and the method has also been used for the modeling of geothermal reservoir stimulation [27]. In addition to its capability for solving multi-physics problems, two features of MPSA allow for the efficient handling of fractures. First, degrees of freedom in MPSA represents cell center displacements; therefore, it does not explicitly imply a spatially continuous approximation. Second, MPSA leads to explicit expressions for traction forces at the grid faces. In the following sections, we will explain how to exploit these features in detail in numerical modeling of deformation of fractured formations. We describe the incorporation of fractures into the MPSA approach for both 2D and 3D spatial discretizations and discuss its convergence properties.

With the goal of developing a numerical method with a broad application area, the MPSA-W method is extended to cover a wide range of mechanical problems in fractured media. We concentrate on three main modeling problems: (1) fracture deformation defined by prescribed displacements on fractures, (2) slip due to the applied tractions at the fracture surfaces and (3) displacements controlled by friction between fracture surfaces. It should be noted that the propagation and initiation mechanisms of fractures are outside the scope of this study.

The study is structured as follows. The model equations and the inclusion of fractures are presented in Section 2. Section 3 starts with the description of the grid structure for the method, followed by the explanation of the numerical discretization approach. We validate the methodology in Section 4 by conducting comparison studies between our implementation and analytical solutions as well as solutions using existing software. In addition, a numerical experiment for deformation of a structure with a complicated fracture network is presented. Finally, we present our conclusions in Section 5.

## 2. GOVERNING EQUATIONS

We consider the intact rock in the subsurface as a linearly elastic medium [1]. However, linear elasticity is not sufficient to model all aspects of deformation for subsurface structures. It is well known that slit-like discontinuities, such as cracks or fractures, are a common type of defect in geological rock [1]. An idealized fracture is described by two surfaces, one on each side of the fracture, subject to specific contact conditions. We model the fractures as two-sided co-dimension one inclusions in the interior of the domain, subject to specific governing equations for their deformation in their tangential and normal directions. Moreover, the deformation of the fractures and the surrounding elastic material are coupled. In the coupling with the deformation of the rock matrix, the surfaces constraining the fractures can be seen as internal boundaries within the elastic domain. Here, we explain our modeling approach in detail, providing a summary of the considered governing equations at the end of the section.

### 2.1. Linear Momentum Balance

The deformation of the intact rock is modeled by static momentum-balance for an elastic medium. For a d-dimensional domain (d=2 or 3), $\Omega$, it is given by:

$$\int_{\partial\Omega} \boldsymbol{T}(\boldsymbol{n})dA + \int_{\Omega} \boldsymbol{f}dV = 0 \text{ on } \Omega, \tag{1}$$
$$\boldsymbol{u} = \boldsymbol{u}^D \text{ on } \partial\Omega^D,$$
$$\boldsymbol{T}(\boldsymbol{n}) = \boldsymbol{T}^N \text{ on } \partial\Omega^N.$$

Here, $\boldsymbol{T}(\boldsymbol{n})$ are the forces on the surfaces of $\Omega$, identified by the outward normal vector $\boldsymbol{n}$; $\boldsymbol{f}$ are body forces acting on the material, and $\boldsymbol{u}$ is the unknown displacement field. The two exponents, $D$ and $N$, denote Dirichlet and Neumann boundary conditions, respectively. In infinitesimal strain theory, the surface forces can be expressed as

$$\boldsymbol{T}(\boldsymbol{n}) = \boldsymbol{\sigma} \cdot \boldsymbol{n}, \tag{2}$$

where $\boldsymbol{\sigma}$ is the Cauchy stress tensor. By introducing $\boldsymbol{\varepsilon}$ as the symmetric part of the deformation gradient $\boldsymbol{\varepsilon} = (\nabla\boldsymbol{u} + (\nabla\boldsymbol{u})^T)/2$, the Cauchy stress tensor can be related to strain through Hooke's Law,

$$\boldsymbol{\sigma} = \mathbb{C}: \boldsymbol{\varepsilon}, \tag{3}$$

where $\mathbb{C}$ is the stiffness tensor. The methods considered herein are developed and implemented for general non-isotropic media as indicated by Eq. (3). However, we will in the specification of the numerical test cases employ the simplification that for isotropic media, the stiffness tensor can be expressed by Lamé parameters $\mu$ and $\lambda$ for; thus, the Cauchy stress tensor can be rewritten in terms of Lamé parameters as

$$\boldsymbol{\sigma} = \mathbb{C}: \boldsymbol{\varepsilon} = 2\mu\boldsymbol{\varepsilon} + \lambda\, tr(\boldsymbol{\varepsilon})\boldsymbol{I}. \tag{4}$$

### 2.2. Fractures as co-dimension one inclusions

We are interested in a problem in which the domain has co-dimension one inclusions, $\Gamma$, that can be considered as discontinuities, i.e., fracture surfaces. Motivated by the method of 'split nodes' presented by [28], the fracture surfaces are considered as line pairs for 2D domains, and pairs of faces for 3D domains, which displace relative to one another, as illustrated in Fig. 1. Following the notation for fracture surfaces from, for example, [12,16], we denote the two sides of the fracture inclusion by subscripts + and -. The tractions on the fracture are defined separately for the positive and negative side. Because of continuity and equilibrium conditions, the relation between traction forces can be written as

$$T_+(n_+) = -T_-(n_-) \text{ on } \Gamma, \qquad (5)$$

where $T_+$ and $T_-$ are the traction forces on the positive and negative side of the inclusion.

In addition to equilibrium conditions between forces on fracture surfaces (Eq. (5)), one more relation for fracture faces is required to complete the system of equations. The presented method can be constructed to include any type of fracture constitutive model for fracture surfaces. In the modeling of subsurface applications, there are three common types of problems one may need to solve: Prescribed displacement jumps over the fracture, prescribed tractions on the fracture, and a frictional relation between normal and tangential forces on the fracture. Here, we focus on and describe these problems in separate subsections.

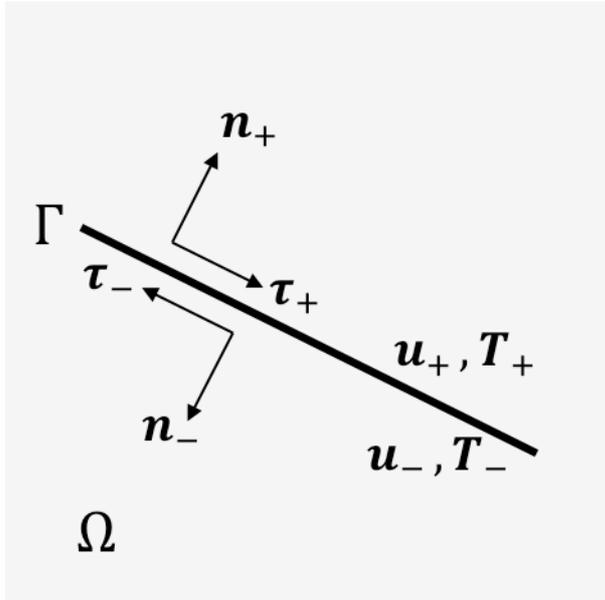

**Fig. 1** Modeling of a fracture. Fractures are modeled as lower dimensional inclusions having two surfaces, labeled + and -.

### 2.2.1. Defined Displacement Jump at the Fracture Surfaces

Slip and dilation, that is, tangential and normal displacement, of a fracture can be modeled by a constitutive model that defines a jump in displacement between positive and – surfaces. A prescribed jump $\Delta u_+^\Gamma$ between the positive and negative sides, which has normal, $\Delta u_{n_+}^\Gamma$, and shear components, $\Delta u_{\tau_+}^\Gamma$, can be written as

$$(u_+ - u_-) = \Delta u_+^\Gamma \text{ on } \Gamma \text{ where } \Delta u_+^\Gamma = n_+ \Delta u_{n_+}^\Gamma + \tau_+ \Delta u_{\tau_+}^\Gamma. \qquad (6)$$

Here, $u_+$ and $u_-$ are the displacements on the positive and negative sides of the fracture surfaces, respectively, and $\tau_+$ denotes the unit vector defining the shear direction of the positive side of the fracture, oriented with respect to $n_+$ according to the right-hand rule. The jump between the fracture faces can be defined according to the negative fracture face, $\Delta u_-^\Gamma$ in similar manner.

### 2.2.2. Defined Traction at the Fracture Surfaces

Next, we consider the case where the tractions at the fracture surface are specified. This setup includes the case of pressurized fracture networks, where fluid pressure acts as a normal force at the fracture surfaces. The defined traction force, $T_+^\Gamma$, at the positive fracture surface can be written as

$$T_+ = T_+^\Gamma \text{ on } \Gamma. \tag{7}$$

The subscript + indicates that the traction is defined on the positive side of the fracture and Eq. (5) provides the traction on the negative side.

### 2.2.3. Fractures Controlled by a Friction Model

In addition to the problem types described in subsections 2.2.1 and 2.2.2, another common type of problem is the friction-controlled fracture displacements in the considered domain. The friction law determines the magnitude of the shear stress on the fracture and is generally motivated by laboratory experiments [29]. The simplest example of a friction model includes constant friction as

$$T_{\tau_+} = \mu_f T_{n_+} \text{ on } \Gamma, \tag{8}$$

using the local coordinate system for a fracture shown in Fig. 1. Here, $T_{\tau_+}$ and $T_{n_+}$ are the magnitudes of the shear and normal tractions on the positive side of the fracture, respectively, and $\mu_f$ is the contact friction between the fracture surfaces.

The friction controlled fracture constitutive model is only meaningful when the fracture surfaces are in contact with each other. Further, even if fracture slip could give rise to dilation of the fracture, we assume, for simplicity, a zero displacement jump in the normal direction when considering this model; that is

$$\boldsymbol{n}_+(\boldsymbol{u}_+ - \boldsymbol{u}_-) = 0 \text{ on } \Gamma. \tag{9}$$

The model requires a definition for the friction coefficient $\mu_f$. Among the numerous proposed models, the most commonly applied in subsurface applications are static-dynamic friction, linear slip-weakening [30], linear time-weakening [31], and Dieterich-Ruina rate-state friction with an aging law [32]. To validate our numerical method, it suffices to consider a constant friction value throughout our simulation; however, the modeling framework applies to more complex friction coefficients as well.

### 2.3. The Complete System of Equations

In summary, we solve the momentum balance equation with the following boundary conditions,

$$\int_{\partial\Omega} \boldsymbol{T}(\boldsymbol{n}) dA + \int_\Omega \boldsymbol{f} dV = 0 \text{ on } \Omega, \tag{10}$$
$$\boldsymbol{u} = \boldsymbol{u}^D \text{ on } \partial\Omega^D,$$
$$\boldsymbol{T}(\boldsymbol{n}) = \boldsymbol{T}^N \text{ on } \partial\Omega^N,$$
$$\boldsymbol{T}_+(\boldsymbol{n}_+) = -\boldsymbol{T}_-(\boldsymbol{n}_-) \text{ on } \Gamma,$$

where the system requires one more condition for the deformation of fracture inclusions. In this study, we consider either of the following conditions:

$$(\boldsymbol{u}_+ - \boldsymbol{u}_-) = \Delta \boldsymbol{u}_+^\Gamma \text{ on } \Gamma, \tag{11}$$
$$\boldsymbol{T}_+ = \boldsymbol{T}_+^\Gamma \text{ on } \Gamma,$$
$$T_{\tau_+} = \mu_f T_{n_+} \text{ on } \Gamma \text{ and } \boldsymbol{n}_+(\boldsymbol{u}_+ - \boldsymbol{u}_-) = 0 \text{ on } \Gamma.$$

While the two first conditions specify a specific displacement jump across the fracture and specific values for the tractions on the different sides of the fracture, the last condition models the shear displacement of the fracture by a friction law.

# 3. THE MPSA DISCRETIZATION WITH FRACTURES

In this section, we first introduce the grid structure that is necessary to discretize the system of equations given in the previous section. The grid structure is the same for all three fracture models summarized in Section 2.3. The section continues with the introduction of MPSA, where we present the discrete form of the momentum-conservation equation along with the discrete forms of the surface stress and the displacement. Moreover, we discuss how to create the discrete system of equations for each type of conditions on the fractures that we consider.

## 3.1. Grid Structure

The simulation domain $\Omega$ is discretized by partitioning it into a set of polyhedral cells that we denote as control volumes $\Omega_i$. We require that the grid conforms to all fracture lines and faces. For cells $i$ and $j$ sharing a boundary, the shared boundary is denoted as face $\omega_{i,j}$. To construct the MPSA discretization, the cells are further divided into one sub-cell per vertex, and all sub-cells associated with a vertex forms an interaction region, see Fig. 2c where the sub-cell associated with vertex $l$ and cell $i$ is denoted $\widetilde{\Omega}_{i,l}$. This will also split faces into sub-faces. We will refer to the sub-faces as $\widetilde{\omega}_{i,j,l}$, where $l$ is the cell vertex associated with the sub-face. In 2D, there will always be two sub-faces per face, whereas in 3D, there will, for example, be three sub-faces per face for triangular faces, and four for rectangular faces.

Among the set of faces of the grid there may be ordinary internal faces, external boundary faces, and internal boundary faces, i.e., fracture faces. To impose fracture relations, such as those defined in Eqs. (5-8), we must represent the displacement on both sides of the fracture. Thus, we create a mesh that includes face pairs corresponding to each side of the fracture. The vertices (edges in 3D) that correspond to tips of the fractures interior to the domain, link the two faces on each side of the fractures. This approach naturally forces the magnitude of the slip to zero at the vertices corresponding to the immersed tips of the fractures. Figure 2 illustrates a conceptual domain that includes fractures and a corresponding mesh for the domain.

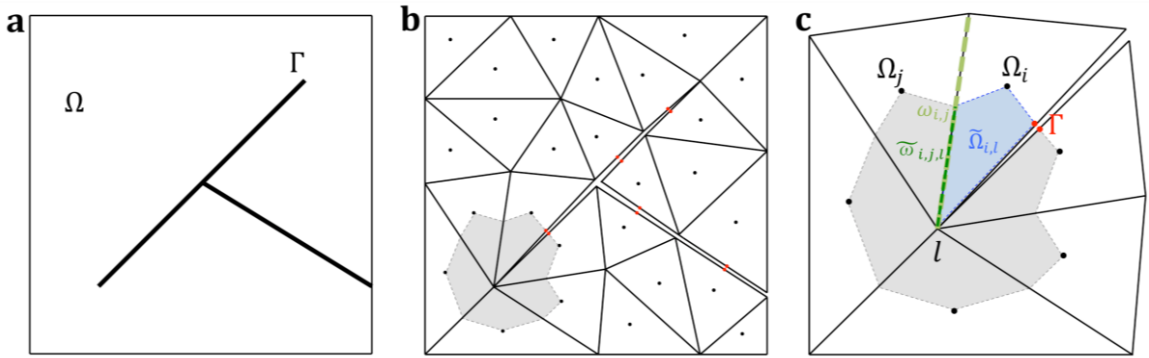

**Fig. 2** The domain and the grid structure illustrated for a 2D case. (a) A 2D conceptual domain including fractures. (b) Conceptual illustration of the grid structure, where the fractures are explicitly represented. Black dots represent the unknown locations of cell centers, while red dots represent the locations of nodes with unknown displacement at fracture faces. For illustration purposes, we show a gap between two red dots, in reality there is no gap between the fracture faces in the grid. Sub-cells creating an interaction region are shown as shaded regions. (c) The details of the grid structure in the interaction region.

## 3.2. Multi-Point Stress Approximation

Considering the grid structure explained above, the momentum conservation for each cell $\Omega_i$ can be written as

$$-\int_{\Omega_i} \boldsymbol{f} dV = \int_{\partial \Omega_i} \boldsymbol{T}(\boldsymbol{n}) dA = \sum_j \int_{\omega_{i,j}} \boldsymbol{T}(\boldsymbol{n}) dA. \quad (12)$$

Furthermore, by introducing the surface stresses between cells $i$ and $j$, $\boldsymbol{T}_{i,j}$, as the stress over face $\omega_{i,j}$, Eq. (12) can be rewritten as

$$-\boldsymbol{f}_i = \frac{1}{|\Omega_i|} \sum_j \boldsymbol{T}_{i,j}, \quad (13)$$

where $-\boldsymbol{f}_i$ is the volume-averaged force over cell $\Omega_i$.

The discretization applies a linear approximation for displacement within each sub-cell $\widetilde{\Omega}_{i,l}$:

$$\boldsymbol{u} \approx \boldsymbol{u}_i + \boldsymbol{G}_{i,l} \cdot (\boldsymbol{x} - \boldsymbol{x}_i), \quad (14)$$

where $\boldsymbol{u}_i$ is the cell-center displacement, $\boldsymbol{x}_i$ is the coordinates of the cell-center, $\boldsymbol{x}$ is a point within the sub-cell $\widetilde{\Omega}_{i,l}$, and $\boldsymbol{G}_{i,l}$ is the gradient in the sub-cell, $\widetilde{\Omega}_{i,l}$. Eq. (14) represents the basis functions created on each sub-cell in terms of two variables: the gradients, $\boldsymbol{G}_{i,l}$, and cell-centre displacements, $\boldsymbol{u}_i$. Combined, these define the discrete displacements in the whole domain. Further, the discrete form of the stresses are written in terms of the gradients, $\boldsymbol{G}_{i,l}$, as

$$\boldsymbol{\sigma}_{i,l} = \frac{\mathbb{C}_{i,l} : \boldsymbol{G}_{i,l} + \left(\mathbb{C}_{i,l} : \boldsymbol{G}_{i,l}\right)^T}{2}, \quad (15)$$

following Hooke's Law.

MPSA methods express the surface stress between cell $i$ and $j$, $\boldsymbol{T}_{i,j}$, as a linear function of the cell-center displacements $\boldsymbol{u}_i$, such that

$$\boldsymbol{T}_{i,j} = \sum_k t_{i,j,k} \boldsymbol{u}_k, \quad (16)$$

where $t_{i,j,k}$ is the stress-weight tensor and $k$ denotes the cells that are neighbors to face $\omega_{i,j}$. More specifically, with each sub-face, $\widetilde{\omega}_{i,j,l}$, we associate a stress weight tensor, denoted as $\hat{t}_{i,j,l,k}$; thus, the stress-weight tensor for face $\omega_{i,j}$ is calculated as the summation of the contributions from each sub-face that is associated with vertex $l$:

$$t_{i,j,k} = \sum_l \hat{t}_{i,j,l,k}. \quad (17)$$

For clarity, the calculation of stress weight tensors (Eqs. (16-17)) and the remaining description of the method is split into two parts: First, a local linear system is created for each vertex of the domain, and gradients are eliminated by solving the local systems. The solution of this local linear system with Eq. (15) leads to calculated stress weight tensors. Then, the global discretization in terms of cell-center displacements is obtained. Several variants of the MPSA approach has been proposed, see [24,22], and, here, we consider the weakly symmetric variant developed in [24]. In the following sections, we start by reviewing this procedure following the weakly symmetric variant of MPSA, before we show how the conditions imposed on the fracture faces affect the local and global systems.

### 3.3. Construction of the Local Linear System without Fractures

As a foundation, we start with the construction of the local systems without fractures. To create a local linear system at each vertex of the domain, two continuity conditions are imposed on each sub-face. First, the continuity of stress over a sub-face is

$$T_{i,j,l} = -T_{j,i,l}, \tag{18}$$

which can be written as

$$\left[\mathbb{C}_i : G_{i,l} + \langle \mathbb{C}_i : G_{i,l} \rangle^T \right] \cdot \bar{n}_{i,j,l} = -\left[\mathbb{C}_j : G_{j,l} + \langle \mathbb{C}_j : G_{j,l} \rangle^T \right] \cdot \bar{n}_{j,i,l}, \tag{19}$$

where $\bar{n}_{i,j,l}$ and $\bar{n}_{j,i,l}$, is the area-weighted normal vector for the corresponding sub-face. The computation of the weak transpose terms, $\langle \mathbb{C}_i : G_{i,l} \rangle^T$ and $\langle \mathbb{C}_j : G_{j,l} \rangle^T$ requires some nuance, see [24].

Similarly, the continuity of displacement over a sub-face can be written as

$$u_i + G_{i,l}(\tilde{x}_{i,m,l} - x_i) = u_j + G_{j,l}(\tilde{x}_{i,m,l} - x_j), \tag{20}$$

where $\tilde{x}_{i,m,l}$ is the calculated continuity point, located one-third of the distance from the face center to the considered vertex for simplex grids, as suggested in [24]. Finally, combining Eq. (19) and (20), the linear system for each vertex in the grid can be expressed as follows:

$$\begin{pmatrix} n^T \mathbb{C} & 0 \\ D_G & D_U \\ 0 & I \end{pmatrix} \begin{pmatrix} G \\ U \end{pmatrix} = \begin{pmatrix} 0 \\ 0 \\ I \end{pmatrix}, \tag{21}$$

where $D_G$ contains distances from the cell centers to continuity points, $D_U$ is a matrix of ±1 representing the contributions from cell centers in Eq. (20), $G$ represents gradients associated with the interaction region, $U$ represents the cell-center displacements, and $I$ is the identity matrix. Observe that the third block row imposes a unit displacement in one component of $U$ at a time. Thus, when solved for the sub-cell gradients $G$, the linear system gives the local deformation response to a unit displacement of each of the cell center variables. The computed sub-cell gradients can be inserted into Eqs. (15-17) to obtain the desired expressions for surface forces in terms of cell-center displacements.

### 3.4. Fracture Implementation in the Local Linear System

We now focus on a local system that includes internal boundary conditions corresponding to fracture faces in the situation where a fracture is in between two regular cells of the domain, see Fig. 3. For the fracture deformation models discussed in Section 3.2, at least one of the continuity equations, Eq. (19) and (20), will have to be replaced. We achieve this by introducing displacement variables at the fracture faces as global variables.

As explained in Section 3.1, the unknowns representing the displacement on the fracture surfaces are located at the face centers. For two neighboring cells $i$ and $s$, split by a fracture, the displacements at the fracture faces can be calculated as

$$u_{i,s,+} = u_i + G_{i,l}(\tilde{x}_{i,m,l} - x_i), \tag{22}$$
$$u_{i,s,-} = u_s + G_{s,l}(\tilde{x}_{s,m,l} - x_s),$$

where $\tilde{x}_{i,m,l}$ and $\tilde{x}_{s,m,l}$ are the split continuity points on sub-face $\tilde{\omega}_{i,m,l}$ and $\tilde{\omega}_{s,m,l}$. If we include Eq. (22) in the former local linear system, Eq. (21), the new local system becomes

$$\begin{pmatrix} \mathbf{n}^T\mathbb{C} & 0 & 0 \\ \mathbf{D}_G & \mathbf{D}_U & 0 \\ \mathbf{D}_G^\Gamma & \mathbf{I} & -\mathbf{I} \\ 0 & \mathbf{I} & 0 \\ 0 & 0 & \mathbf{I} \end{pmatrix} \begin{pmatrix} \mathbf{G} \\ \mathbf{U} \\ \mathbf{U}^\Gamma \end{pmatrix} = \begin{pmatrix} \mathbf{0} \\ \mathbf{0} \\ \mathbf{0} \\ \mathbf{I} \\ \mathbf{I} \end{pmatrix}, \qquad (23)$$

where $\mathbf{D}_G^\Gamma$ contains distances from the centers of the cells that have internal boundaries to the continuity points. Here, the third row represents the displacements on the internal boundaries corresponding to fracture faces (Eq. (22)) and $\mathbf{U}^\Gamma$ represents both of the displacements at the positive and negative fracture faces. Again, the two last rows of Eq. (23) impose a unit displacement in each cell-center and fracture-face at a time. By solving Eq. (23) for the sub-cell gradients $\mathbf{G}$, we can obtain the displacement gradients in terms of both cell-center displacements and the displacements defined at the fracture faces. The gradients can be inserted into Eqs. (15-17) to obtain the desired expressions for surface forces in terms of both cell-center displacements and displacements on the internal boundaries corresponding to the fracture.

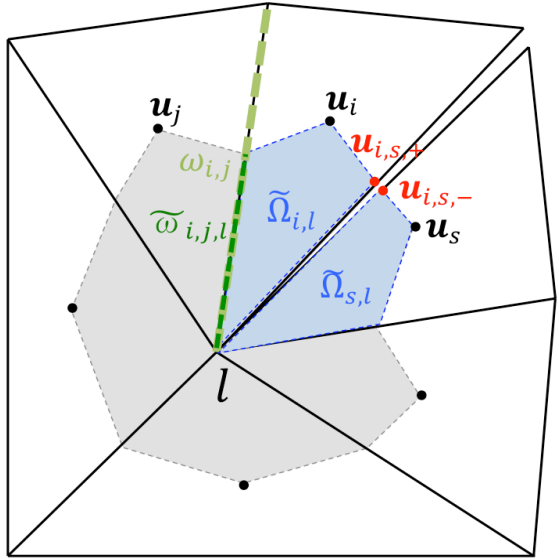

**Fig. 3** Fracture implementation in the local linear system. The fracture faces are located between cells $i$ and $s$. While the unknowns $\mathbf{u}_i$ and $\mathbf{u}_s$ represent the displacements in the cells centers of $i$ and $s$, the unknowns $\mathbf{u}_{i,s,+}$ and $\mathbf{u}_{i,s,-}$ represent the discrete displacements on the fracture surfaces. For clarity, the colors of light green, dark green, and blue are associated with the faces, sub-faces, sub-cells, in an example interaction region, respectively.

### 3.5. Construction of the Global Linear System

By solving the local systems for each vertex in the grid, the discrete tractions can be obtained in the desired form of Eq. (16). Furthermore, a discrete version of the momentum balance equation, Eq. (1), for each cell $\Omega_i$ can be written accordingly. By combining the balance equations for all cells, the final form of the force-balance equation in terms of displacements can be written as

$$-|\Omega_i|\mathbf{f}_i = \sum_j \mathbf{T}_{i,j} = \sum_j \sum_k t_{i,j,k} \mathbf{u}_k. \qquad (24)$$

For a domain without fracture inclusions, one can write the linear system of Eq. (24) in the form of

$$\psi^\Omega U = 0, \qquad (25)$$

where the matrix $\psi^\Omega$ represents the terms of $\sum_j \sum_k t_{i,j,k}$ in Eq. (24).

In the presence of fracture inclusion, Eq. (25) is augmented according to the fracture constitutive models that are discussed in Section 2.2. To include the continuity of the stresses (Eq. (5)) between each side of the fracture faces (+, -) to the global linear system, Eq. 25, we make use of the definition of the traction forces in terms of cell center displacements (Eqs. (15-17)). The stress weight tensors of one particular face have contributions from all neighboring sub-cells (see Fig. 3); thus we represent stress weight tensors on positive fracture faces (($\sum_k t_{i,s,+,k}$)) in the form of $\tilde{\gamma} = (\tilde{\gamma}^\Omega \quad \tilde{\gamma}^\Gamma_+ \quad \tilde{\gamma}^\Gamma_-)$ and the traction vectors for positive fracture faces can be calculated as

$$(\tilde{\gamma}^\Omega \quad \tilde{\gamma}^\Gamma_+ \quad \tilde{\gamma}^\Gamma_-) \begin{pmatrix} U \\ U^\Gamma_+ \\ U^\Gamma_- \end{pmatrix}. \qquad (26)$$

The equilibrium condition between the positive and negative side of the fracture can be written as,

$$(\tilde{t}^\Omega \quad \tilde{t}^\Gamma_+ \quad \tilde{t}^\Gamma_-) \begin{pmatrix} U \\ U^\Gamma_+ \\ U^\Gamma_- \end{pmatrix} = \mathbf{0}, \qquad (27)$$

where $\tilde{t} = (\tilde{t}^\Omega \quad \tilde{t}^\Gamma_+ \quad \tilde{t}^\Gamma_-)$ represents the summation of the stress weight tensors between positive and negative fracture faces (($\sum_k t_{i,s,+,k} + \sum_k t_{i,s,-,k}$)). Now we are equipped with combine the global linear system, Eq. (24), with the continuity of the stresses (Eq. (5)) between each side of the fracture faces (+, -). Further, we introduce matrices $L$ and $R$, with $L = (L^\Omega \quad L_+ \quad L_-)$ decomposed into the contribution from the interior, and the two sides of the fractures. Eq. (25) is augmented to the form of

$$\begin{pmatrix} \psi^\Omega & \psi^\Gamma_+ & \psi^\Gamma_- \\ \tilde{t}^\Omega & \tilde{t}^\Gamma_+ & \tilde{t}^\Gamma_- \\ L^\Omega & L_+ & L_- \end{pmatrix} \begin{pmatrix} U \\ U^\Gamma_+ \\ U^\Gamma_- \end{pmatrix} = \begin{pmatrix} 0 \\ 0 \\ R \end{pmatrix} \qquad (28)$$

for a system that includes fractures. Here, the terms with superscript Ω are associated with the cell centers, while superscript Γ represents the terms associated with fracture faces. Likewise, the subscripts + and – denote terms associated with the positive and negative sides of a fracture-face pair, respectively. The first row in Eq. (28) is the force-balance equations for each cell in the grid, namely, Eq. (25). The second row represents the continuity of the stresses (Eq. (5)) between each side of the fracture. The $L$ and $R$ matrices in the third row are constructed according to the relevant model for fracture deformation. In the following, we will show how to create $L$ and $R$ matrices for the three most common problems: defined displacement jump between fracture faces (Eq. (6)), prescribed traction at the fracture faces (Eq. (7)), and modeling of friction-controlled fractures (Eqs. (8-9)).

### 3.5.1. Defined Displacement Jump at the Fracture Surfaces

In this case, the fracture surfaces are modeled such that the surfaces are displaced relative to each other by a constant amount, $\Delta u^\Gamma_{i,s}$, as in Eq. (6). The $L$ and $R$ matrices in Eq. (28) should be written as

$$L = (\mathbf{0} \quad I \quad -I), R = \Delta u^\Gamma_{+_{i,s}}. \qquad (29)$$

As we see from Eq. (29), when $\Delta u^\Gamma_{+_{i,s}}$ is equal to $\mathbf{0}$, the system disregards the existence of fractures.

### 3.5.2. Defined Traction at the Fracture Surfaces

Similarly, one would like to approximate the displacement distribution in the domain when the tractions are applied to the fracture faces. For the traction, $T^{\Gamma}_{+_{i,s}}$, is defined to the positive side of the fracture surfaces, the solution approximation of this problem can be accomplished by creating the $L$ and $R$ matrices such that

$$L = (\widetilde{\boldsymbol{\gamma}}^{\Omega} \quad \widetilde{\boldsymbol{\gamma}}^{\Gamma}_{+} \quad \widetilde{\boldsymbol{\gamma}}^{\Gamma}_{-}), R = T^{\Gamma}_{+_{i,s}}. \tag{30}$$

### 3.5.3. Fractures Controlled by a Friction Model

In addition to the problems discretized previously (with defined jump or traction at the fracture surfaces), a constitutive friction model may control the tractions of the fracture surfaces. In this problem, it is assumed that the shear stress is controlled through Eq. (8). If we assume that a constant friction, $\mu_f$, controls the deformation behavior of the fracture, the corresponding forms of the $L$ and $R$ matrices become

$$L = \begin{pmatrix} |\widetilde{\boldsymbol{\gamma}}^{\Omega}_{\tau_+}| - \mu_f |\widetilde{\boldsymbol{\gamma}}^{\Omega}_{n_+}| & |\widetilde{\boldsymbol{\gamma}}^{\Gamma}_{\tau_+}| - \mu_f |\widetilde{\boldsymbol{\gamma}}^{\Gamma}_{n_+,+}| & |\widetilde{\boldsymbol{\gamma}}^{\Gamma}_{\tau_+,-}| - \mu_f |\widetilde{\boldsymbol{\gamma}}^{\Gamma}_{n_+,-}| \\ 0 & I & -I \end{pmatrix}, R = \begin{pmatrix} 0 \\ 0 \end{pmatrix}. \tag{31}$$

The subscripts of $n_+$ and $\tau_+$ represents the normal and shear components of the $\widetilde{\boldsymbol{\gamma}}$. This system is different from the previous sets of equations, in that the first row in $L$ includes nonlinear relations constructed by Eq. (8), and the second row ensures zero normal displacements across the fracture surfaces according to Eq. (9).

In this study, all above-mentioned systems of linear equations are solved with a direct solver as implemented by Matlab's backslash operator. For the nonlinear system, Newton's method [33] is implemented.

## 4. NUMERICAL SIMULATIONS

The performance of the method is verified with several comparison studies for 2D and 3D domains. We start by examining the convergence properties of the method with a 2D domain that contains only one fracture, located approximately at the middle. Moreover, an example that includes a complicated fracture distribution is also presented to show the capabilities of the method. Finally, the convergence of the method for a 3D case is examined. In this case, we analyze the convergence of the method relative to an approximate reference solution, which is found by a finite-element discretization of the problem. For simplicity, both Lamé parameters of the medium are set to 1 in all numerical simulations.

The convergence of the proposed method is examined by calculating the discrete L2-norm of the displacement,

$$err_{\boldsymbol{u}_i} = \frac{\left(\sum_{\Omega_i} m_{\Omega_i} (\boldsymbol{u}_{ref} - \boldsymbol{u}_i)^2\right)^{\frac{1}{2}}}{\left(\sum_{\Omega_i} m_{\Omega_i} (\boldsymbol{u}_{ref})^2\right)^{\frac{1}{2}}}, \tag{32}$$

where $m_{\Omega_i}$ is the area of the considered cell and $\boldsymbol{u}_{ref}$ is the reference solution. The traction error is also calculated in a similar manner as

$$err_{\boldsymbol{T}_{i,j}} = \frac{\left(\sum_{\omega_{i,j}} m_{\omega_{i,j}} (\boldsymbol{T}_{ref} - \boldsymbol{T}_{i,j})^2\right)^{\frac{1}{2}}}{\left(\sum_{\omega_{i,j}} m_{\mathcal{F}_{\omega_{i,j}}} (\boldsymbol{T}_{ref})^2\right)^{\frac{1}{2}}}, \tag{33}$$

where $m_{\omega_{i,j}}$ is the area of the considered face and $\boldsymbol{T}_{ref}$ is the reference solution.

## 4.1. Validation for Single-Fracture Case

In this part, the validation examples for the three types of problems will be addressed: defined displacement jump at the fracture surfaces (case 1), defined traction at the fracture surfaces (case 2), and fractures controlled by a friction model (case 3). The analytical solutions for each type of problem are available in the literature. However, the analytical solutions are defined in an infinite domain. Therefore, there are limitations in using analytical solutions as reference solutions to examine the convergence properties of the method for the cases in which the boundary effect is significant. In these cases, we use the approximate solution with a very fine grid as the reference solution. The convergence studies for all problems are conducted with the same domain. The total length of the domain is 50 m × 50 m with a fracture located at the center the domain. The fracture is 10 m in length and is at a 20° angle from the x-axis. The domain is discretized using a non-structured simplex grid that is created by constrained Delaunay triangulations [34]. On the left side of Fig. 4, we show an example of the coarsest grid, which has four face pairs to discretize the fracture surfaces. For convergence studies, refinement of the grid is conducted such that the numbers of face pairs for the fracture discretization are 8, 16, 32, 64, and 128. The total number of cells in the domain ranges from approximately 1,500 for the coarsest grid to 1.5 million for the finest. To minimize the dependence on a particular unstructured grid sequence, 10 independent grids are constructed, keeping the number of face pairs on the fracture surfaces constant. Then, the evaluation is extended by calculating the average error for each group of grids with the same number of face pairs. We also use non-hierarchical grid refinement for our analysis. Figure 4 (right) shows examples of the non-hierarchical grid refinement at the same tip area for the first three grid refinements. The refinement is performed depending only on the maximum area of each cell; each time the grid is refined, a new grid is created. We will comment on the effect of non-hierarchical refinement in the results.

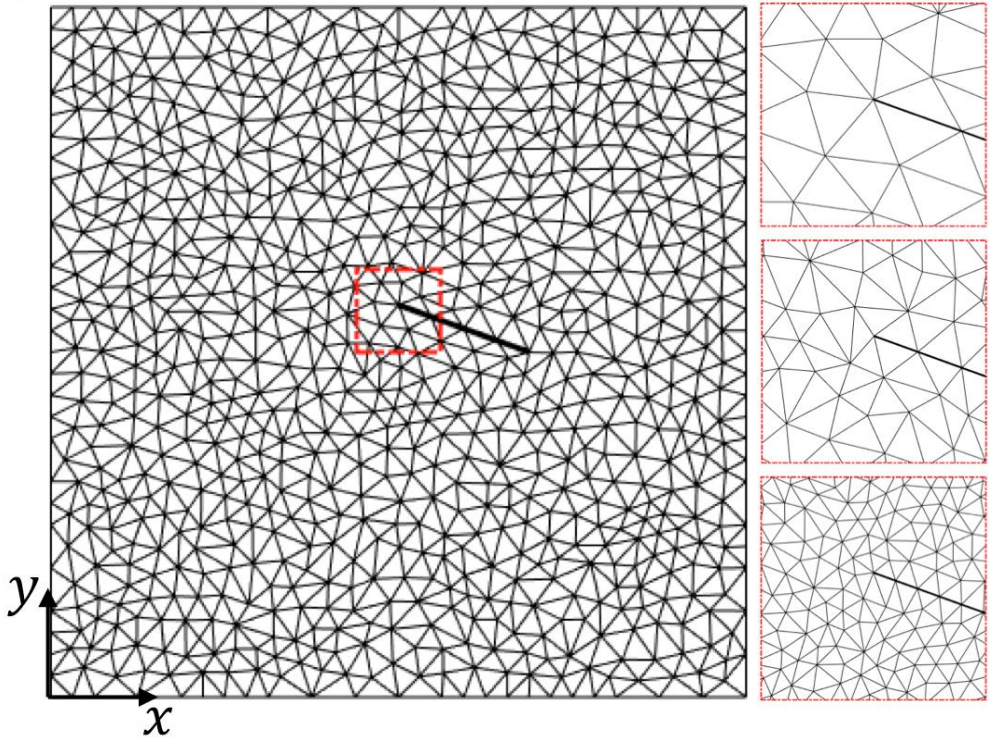

**Fig. 4** The example grid and the non-hierarchical grid refinement that are used for the problems, which have analytical solutions for single cracks. The thin lines show the gridding, while the thick ones represent the fracture. (Left) The coarsest grid used for the validation of the method. (Right) The non-hierarchical grid refinement examples shown close to the tip of the fracture. (The region bounded with the red dotted line on the left figure). Each time the grid resolution is increased, a new grid is created.

### 4.1.1. Case 1: Defined Displacement Jump at the Fracture Surfaces

The first validation problem consists of an infinite, two-dimensional, homogeneous, isotropic, elastic medium that has a constant displacement discontinuity over a finite line segment. It is assumed that the displacements are continuous everywhere except the line segment. The line segment is considered a fracture that has positive and negative sides as described in Section 2.2. The constant displacement discontinuity in the tangential direction is defined between the positive and negative sides of the fracture by setting the tangential component of $\Delta \boldsymbol{u}_+^\Gamma$ in Eq. (6) to 0.001. A detailed description of the problem and the derivation of the analytical solution can be found in [12,13]. The induced displacements and stresses at any point $(x, y)$ in the domain caused by the constant discontinuity are given in the appendix. The traction vectors are calculated for each face using stress values with the face normal vectors.

Figure 5 shows the displacement errors calculated in the whole domain for each grid and the average error. The method is approximately first-order convergent in displacement for this problem, which is lower than in previous numerical studies of MPSA methods without fractures. The reduced convergence order is expected when the domain includes discontinuity. Moreover, since the averages of the displacement errors do not deviate dramatically from the errors calculated for each grid, we can also conclude that the grid structure has only a minor effect on the displacement convergence of the problem.

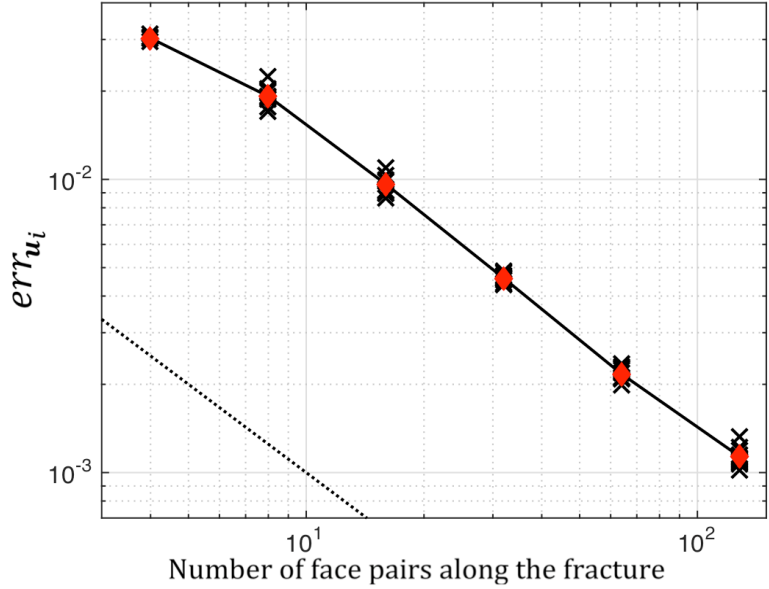

**Fig. 5** The convergence behavior for displacement in the whole domain for case 1. The red diamond shows the average of the error for each group of grids with the same number of face pairs. The dotted line is the 1st-order reference line.

Figure 6 shows the convergence behavior for tractions in the whole domain. As in the previous figure, the line with diamonds in Fig. 6 shows the average values. As we increase the number of DOFs in the domain, the stress does not converge. This behavior is common for lowest order numerical approximations, which fail to capture the singularity at the fracture tips that are inherent in linear elastic fracture mechanics [35,36]. When the number of DOFs is increased, the locations of the unknowns approach the singularity at the tip, and the error becomes more dominant.

To measure the accuracy of the stress approximation, we quantify the stress convergence by excluding values within an estimated inelastic zone, where the linear elastic approximation cannot be expected to hold true, see e.g. [37], [7] for details. The radius of the inelastic zone in general depends on the yield strength of the material. Herein, the tip zone radius is simply set to 0.12, that is, much smaller than the fracture length. When values inside this radius are excluded from the stress calculation, convergence is recovered. This effect is shown in Fig. 6 by the line with stars. In addition, since the averages of the stress errors do not deviate significantly from the errors calculated for each grid, we can again conclude that the grid structure has only a minor effect on the stress convergence of the problem.

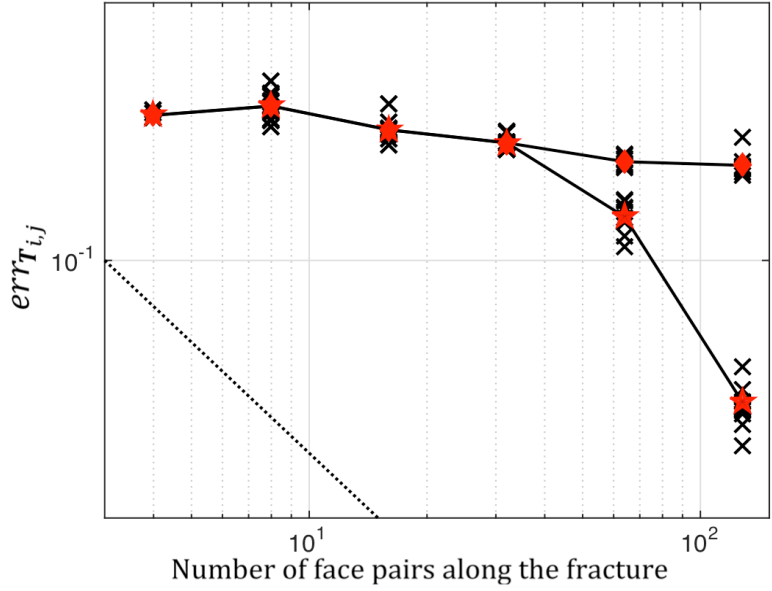

**Fig. 6** The convergence behavior for stress in the whole domain for case 1. The line with diamonds shows the average stress convergence behavior dominated by error caused at the fracture tips, whereas the line with stars shows the stress convergence behavior when the tips of the fracture are eliminated. The dotted line is the 1st-order reference line.

### 4.1.2. Case 2: Defined Traction at the Fracture Surfaces

The second validation problem is defined in the same domain as case 1. This time, the fracture is subject to a constant pressure along its surfaces. The constant pressure can be considered as constant traction in the direction normal to the fracture faces. An analytical solution for the opening between the fracture surfaces (i.e., the normal relative displacement between the fracture surfaces) was derived by Sneddon [38] as

$$\Delta u_{n_+}^{\Gamma}(\eta) = \frac{(1-v)P\Delta L}{\mu}\sqrt{1 - \frac{\eta^2}{(\Delta L/2)^2}}, \qquad (34)$$

where $v$ is Poisson's ratio, $\mu$ is the shear modulus, $P$=0.001 is the applied pressure, $\Delta L$ is the fracture length, and $\eta$ denotes distance from the centre of the fracture, $-\Delta L/2 \leq \eta \leq \Delta L/2$. Note that the traction forces are applied to the fracture surfaces only in the normal direction; therefore, while Eq. (34) gives the normal component of the displacement jump, the shear components of the displacement jump are not defined. The comparison between the analytical solution and the solution approximated using 16 face pairs is shown in Fig. 7. The approximate solution is consistent with the analytical solution.

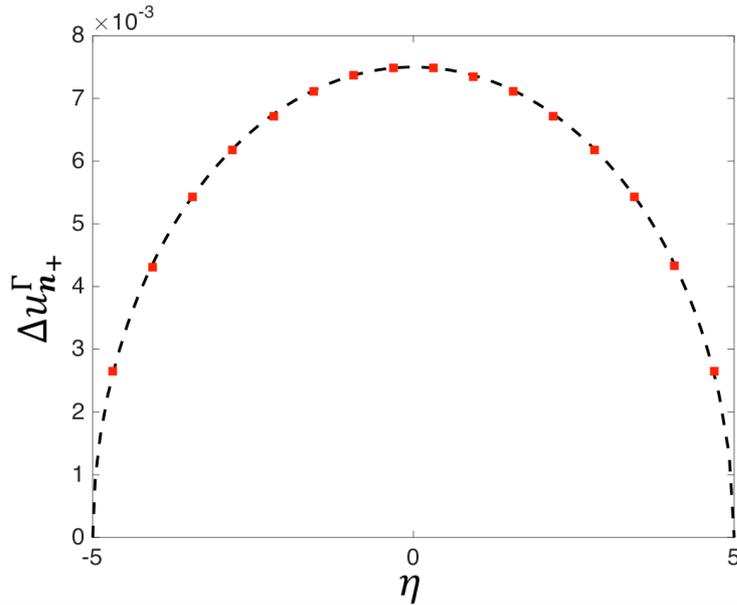

**Fig. 7** Comparison between the analytical solution and the solution approximated using 16 face pairs for the pressurized fracture problem (case 2). The dashed line shows the analytical solution, while the red dots show the numerical solution.

Although the approximate solution is consistent with the numerical solution, the grid structure and the singularity of the stress at the fracture tip have dominant effects on the solution since the forces are defined on the fracture faces. This effect can be clearly observed in Fig. 8, where both the error of all grids and their average is shown. Note that we observed smaller deviations from the average value when we used structured grids at the fracture tips. The method provides 1st-order convergence on average in this case as well.

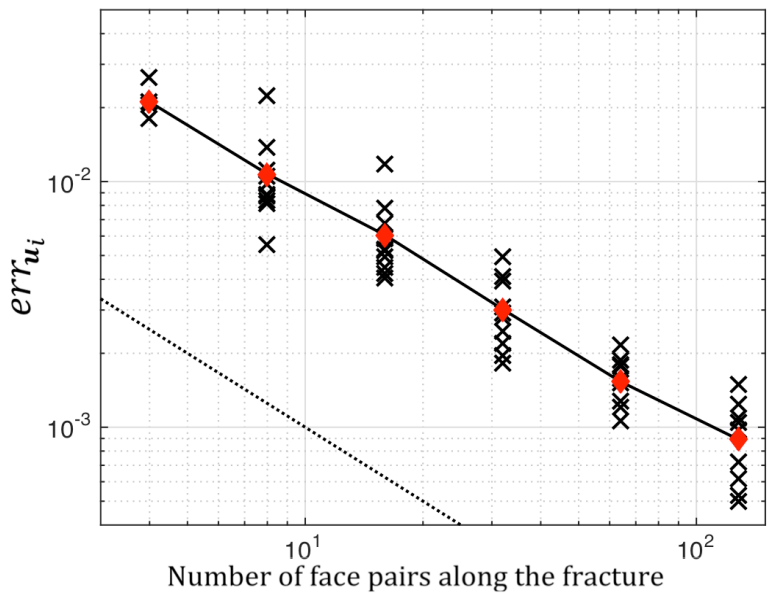

**Fig. 8** The average convergence behavior for displacement at the fracture for case 2. The red diamond shows the average of the error for each group of grids with the same number of face pairs. The dotted line is the 1st-order reference line.

### 4.1.3. Case 3: Fractures Controlled by a Friction Model

In this case, the frictional sliding of a single fracture is considered. We use the same computational domain as in cases 1 and 2. The angles of the fracture inclination and the friction are 20° and 30°, respectively. For the current problem, we wish to create boundary conditions such that the fracture is subject to a compressive stress of $\sigma_{xx}$, while the physical contact between fracture faces is maintained during the simulation. The considered stress conditions can be written as

$$\sigma_{xx} = 0.001, \quad (35)$$
$$\sigma_{yy} = 0,$$
$$\sigma_{xy} = 0.$$

Figure 9 illustrates the boundary conditions for this case. We apply these stress conditions as Neumann boundary conditions on the boundary $\partial\Omega^N$. To avoid an indefinite problem, we define a Dirichlet boundary condition on one side of the domain, $\partial\Omega^D$. For $\partial\Omega^D$, the displacement at the origin and in the $x$-direction, $u_x$, is set to zero, and the displacement in the $y$-direction, $u_y$, is found with the help of Hooke's law, 2D strain-displacement relations, and Eq. (35) as

$$\sigma_{yy} = \lambda\varepsilon_{xx} + (\lambda + 2\mu)\varepsilon_{yy} \text{ where } \varepsilon_{xx} = \frac{\partial u_x}{\partial x}, \varepsilon_{yy} = \frac{\partial u_y}{\partial y}, \quad (36)$$

$$\frac{\partial u_x}{\partial x} = -\frac{(\lambda + 2\mu)}{\lambda}\frac{\partial u_y}{\partial y},$$

$$\sigma_{xx} = \lambda\varepsilon_{yy} + (\lambda + 2\mu)\varepsilon_{xx},$$

$$\frac{\partial u_y}{\partial y} = -\frac{0.001\lambda}{4\mu(\lambda + \mu)}.$$

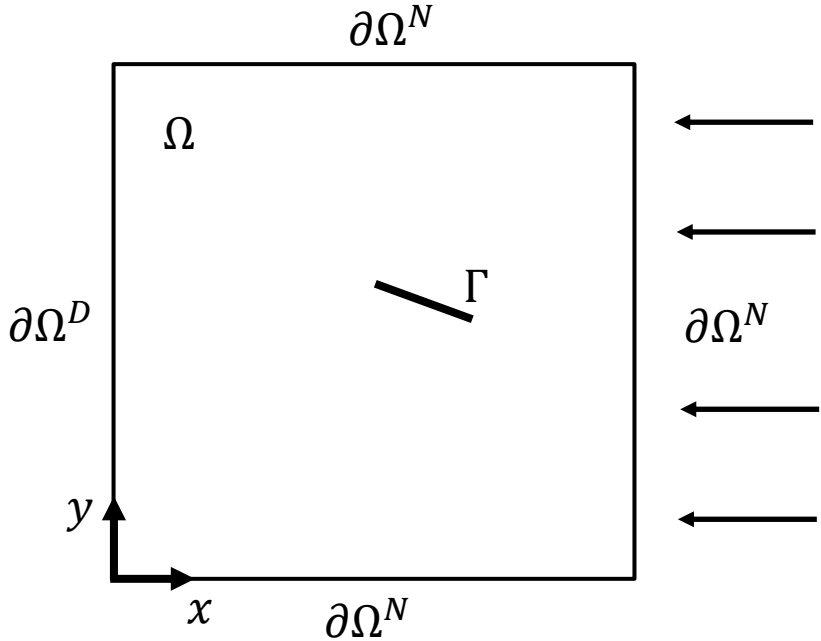

**Fig. 9** Boundary conditions used for the case in which the fracture deformation is controlled by a friction model. The fracture is under compressive stress.

In this problem, the convergence of the numerical solution is examined using the solution with the finest grid as the reference solution. The convergence rate for the displacement is shown in Fig. 10. As in the previous cases, the average convergence is of 1st order. Moreover, the averages of the displacement errors exhibit only minor deviation from the errors calculated for each grid.

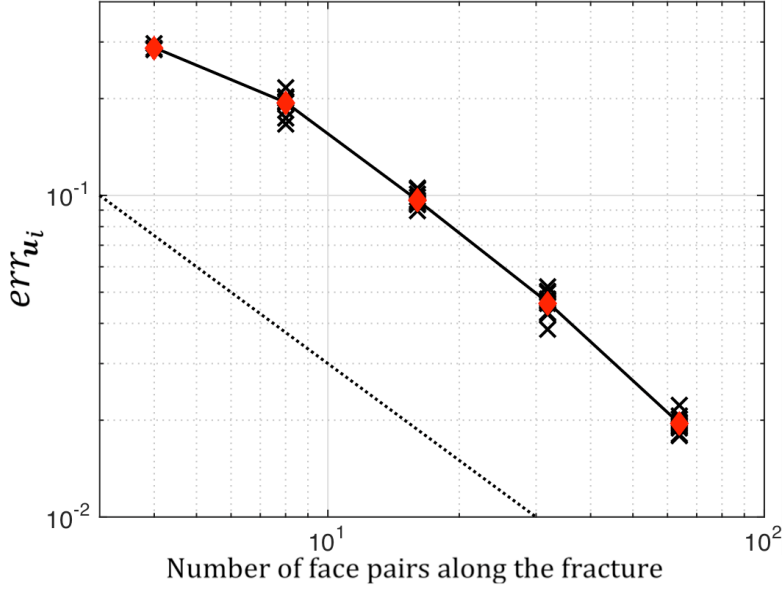

**Fig. 10** The convergence behavior for displacement at the fracture for the nonlinear problem. The red diamond shows the average of the error for each group of grids with the same number of face pairs. The dotted line is the 1st-order reference line.

In addition to the convergence analysis shown in Fig 10, we also make a qualitative comparison between an existing analytical solution for the problem defined in an unbounded domain and the approximate solution defined in a bounded domain. Figure 11 illustrates the geometrical configuration of the analytical problem. The analytical solutions for the shear and normal tractions and the shear component of the displacement jump between the fracture surfaces are given in [39] as

$$T_{n_+} = -\sigma_{xx} \sin^2 \alpha, \tag{37}$$

$$\Delta u^{\Gamma}_{\tau_+}(\eta) = \frac{4(1-\nu^2)T_{\tau_+}}{E} \sqrt{\left(\frac{\Delta L}{2}\right)^2 - \left(\eta - \frac{\Delta L}{2}\right)^2},$$

$$T_{\tau_+} = \sigma_{xx} \sin \alpha \, (\cos \alpha - \sin \alpha \tan \phi),$$

where $0 \leq \eta \leq \Delta L$, $T_{n_+}$ is the normal and $T_{\tau_+}$ is the shear traction on positive face, $\alpha$ is the fracture inclination angle, $\phi$ is the friction angle, and $E$ is Young's modulus.

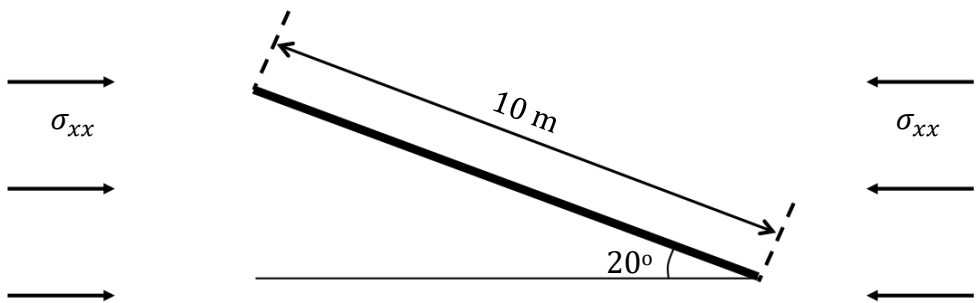

**Fig. 11** A crack under compression in an unbounded domain.

Figure 12 shows the analytical solution and the numerical solution of the shear displacement between the fracture surfaces. As expected, there is only a small difference between the analytical and numerical solutions. The main reason for the difference between the two solutions is the discrepancy between boundary conditions. The number of iterations of the nonlinear solver is 2 for this problem.

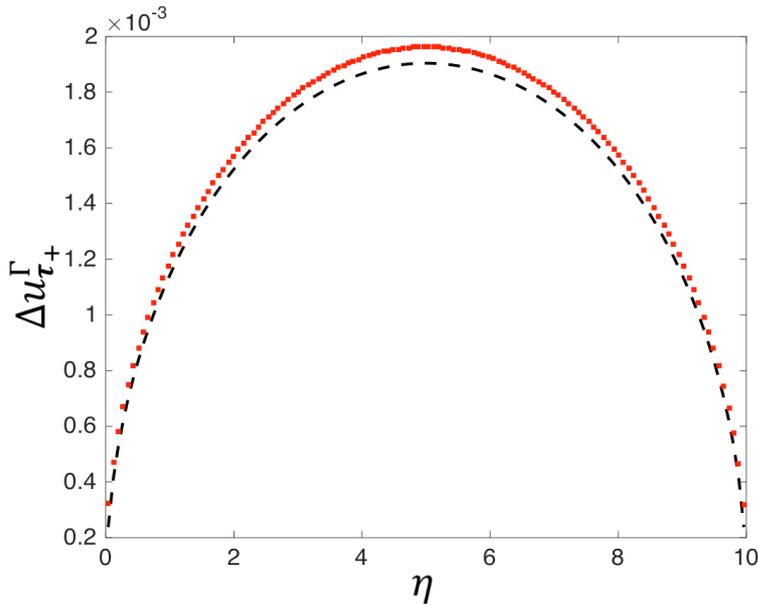

**Fig. 12** The shear displacement controlled by friction. The analytical solution is shown as a dashed line, whereas the numerical solution is shown as red dots. For this numerical solution, 124 face pairs are used.

### 4.2. Complex Numerical Examples

We now consider more complex fracture geometries. We first examine the convergence properties of the method when the domain includes several complications in the numerical sense. Then, we examine a 3D domain with a fracture in the middle. We also show the convergence behavior for the 3D example, based on a reference solution that is approximated using the finite-element method.

#### 4.2.1. Numerical Example with Intersecting Fractures

After the validation of the method illustrated in Section 4.1, a more challenging domain, shown in Fig. 13, is examined to demonstrate the capabilities of the method. The domain is discretized using non-structured simplex grids that are created by constrained Delaunay triangulations [40]. The domain includes challenges such as fracture intersections of several types, a fracture passing through the domain, and a disconnected fracture network. These numerical challenges are handled straightforwardly by the model as the variables for the fracture faces are defined at the centers of the grid faces, as described in Section 3. In addition, the method can accommodate fractures that have immersed tips and fractures that cross the boundary, if the grid is created accordingly (see Section 3.1). Here, we prescribe a constant displacement jump between fracture faces. The displacement is defined in the shear direction of each fracture. Figure 14 shows the deformation results for two grids with 47 and 160 face pairs to discretize the fractures. The colors indicate the magnitudes of the deformations, and the directions of the deformations are shown with arrows. The sizes of the arrows are also related to the deformation magnitudes. The results from using the coarse and fine meshes are consistent.

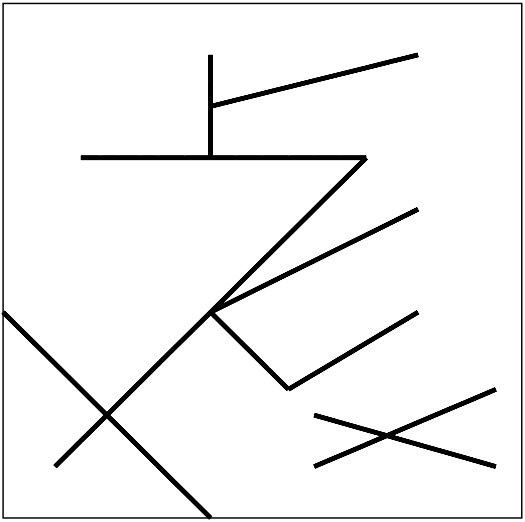

**Fig. 13** The domain includes several intersections as intersecting fractures and disconnected fracture networks.

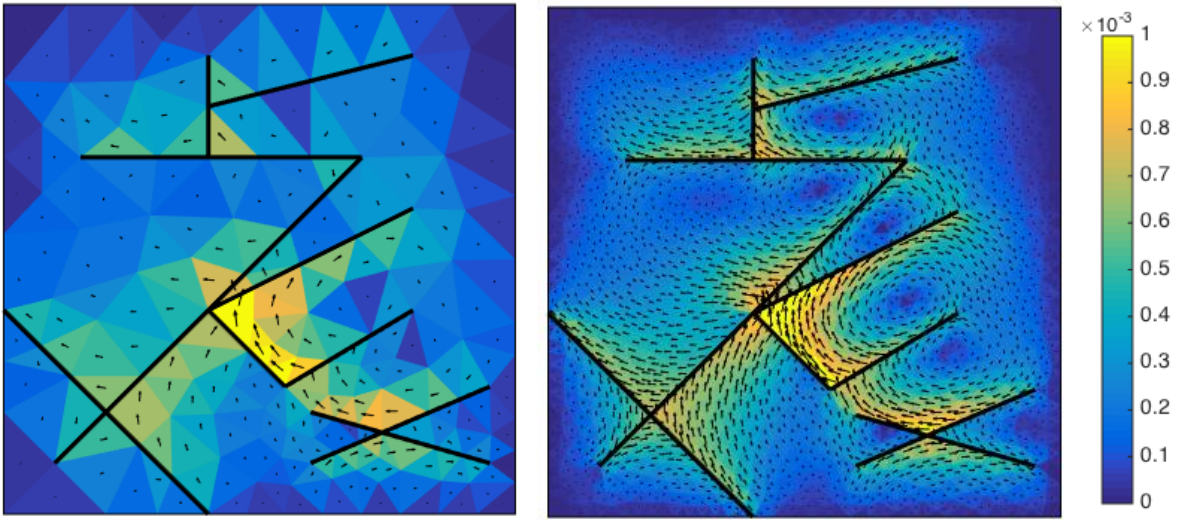

**Fig. 14** Approximate deformation distributions for the solution using 47 face pairs to discretize the fractures (left) and the solution using 160 face pairs to discretize the fractures (right).

To conduct a convergence examination for this example, five different grids that have 47, 83, 160, 325, 671 numbers of face-pairs are used. The reference solution is obtained using a relatively fine grid that has approximately 1,300 face pairs to discretize the fractures. The convergence study is performed by comparing the tractions at the fracture faces between solutions obtained using coarser grids and the reference solution. Figure 15 shows the convergence ratio of the problem. The convergence ratio for this problem is higher than 1st order.

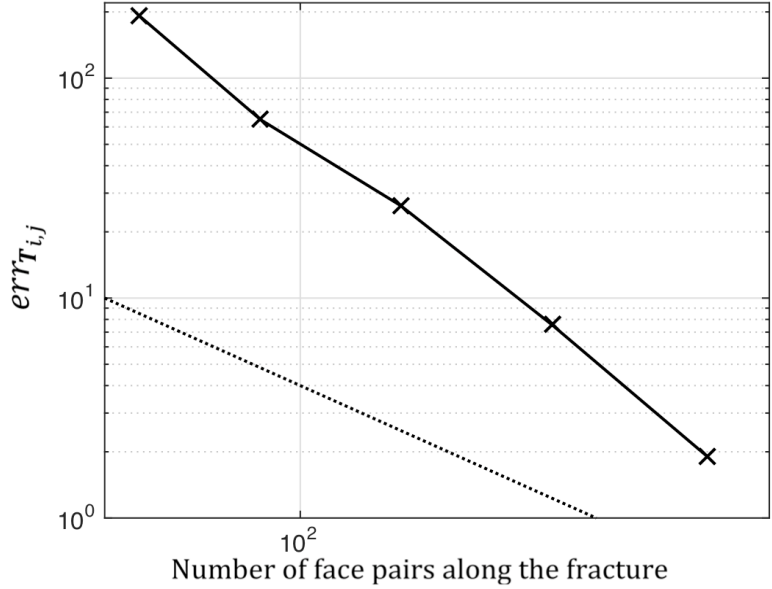

**Fig. 15** The convergence behavior for the problem includes several numerical challenges. The dotted line is the 1st-order reference line.

### 4.2.2. Comparison Example with Finite-Element Method (3D)

As a final example, we examine a 3D domain with a fracture in the middle. A lateral shear displacement jump at the fracture surfaces is defined. The domain and the defined jump are shown in Fig. 16. The displacement jump is defined such that its maximum value is attained at the middle of the fracture and is zero at the edges of the fracture. The equation defining the shear component of the jump can be written as

$$\Delta u_{\tau_+}^{\Gamma}(\eta_x, \eta_y) = 0.001 \left( \sin\left(\frac{\pi \eta_x}{\Delta L}\right) \sin\left(\frac{\pi \eta_y}{\Delta L}\right) \right) \tag{38}$$

where $0 \leq \eta_x \leq \Delta L$ and $0 \leq \eta_y \leq \Delta L$ and the length of each edge of the fracture, $\Delta L$, is 36. The local coordinate system relative to the fracture is shown in Fig. 16.

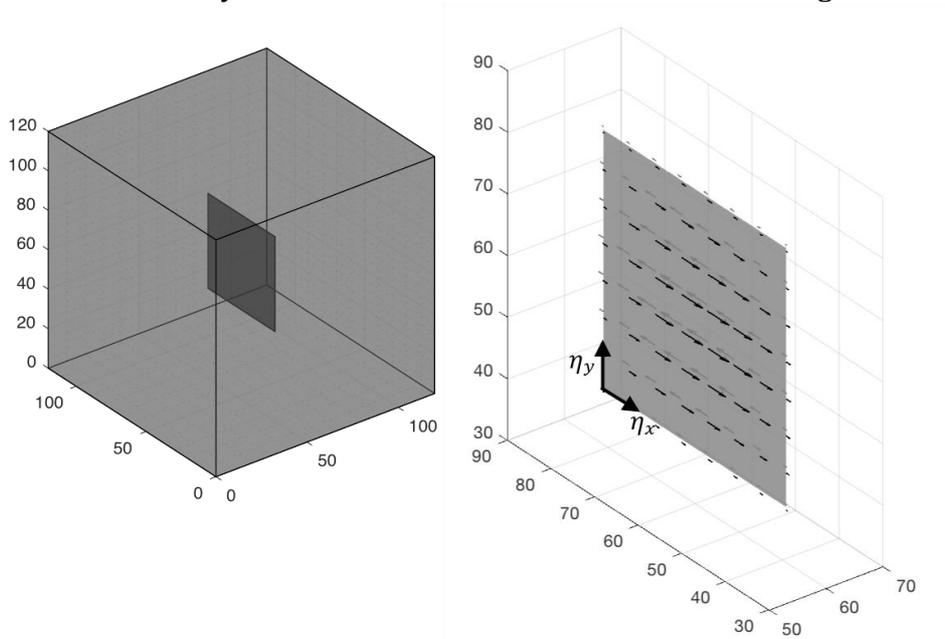

**Fig. 16** The considered 3D domain with a fracture (left). The fracture is shown in dark grey in the middle. Distribution of the prescribed jump at the fracture faces (right).

For this case, we compare the convergence of MPSA and a FEM by using Cartesian hexahedron grids. The number of cells in the domain ranges from 1000 to 216000 for both MPSA and FEM. The reference solution for this problem is assumed to be the solution calculated with FEM by using the finest grid (216000 cells). We use the open-source software package Pylith [16,41], which has been developed for the simulation of crustal deformations. Similar to the present method, Pylith uses cohesive cells and adds DOFs to the grid to define the relative motion of the fracture surfaces. Figure 17 shows the convergence for displacement for both the current method and FEM (Pylith). The convergence ratio is slightly higher than 1 for FEM (Pylith) and slightly lower than 1 for MPSA, while the errors are of the same magnitude for both methods. The improved convergence order of the finite-element method may be related to the reference solution being calculated with this method.

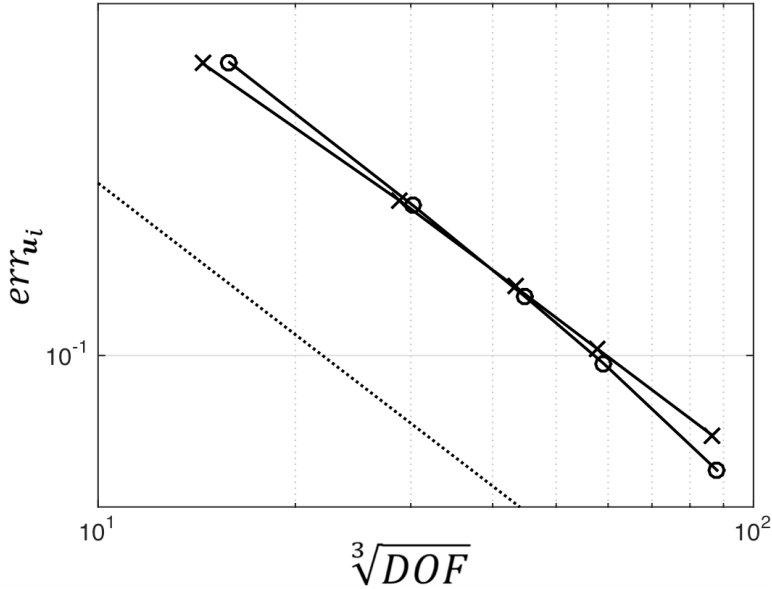

**Fig. 17** Convergence behavior for displacement in the 3D problem. The line with circles shows the Pylith (FEM) results, and the line with crosses shows the MPSA results. The dotted line is the 1st-order reference line.

## 5. CONCLUSION

Motivated by the modeling of subsurface applications, an FVM for the deformation of fractured media, in which pre-existing fractures are represented explicitly in the model, is proposed. The method is constructed based on MPSA with the same grid structure as the FVM for flow, e.g., multi-point flow approximations. The method is therefore particularly well suited for problems involving coupled flow and deformation in fractured formations.

The fractures in a medium physically represent internal discontinuities, and the physics of processes related to fractures can differ and can be modeled by different governing equations than those of the rest of the medium. We consider three models for the deformation of the fractures coupled to the rest of the medium: (1) defined displacement jump at the fracture surfaces, (2) defined traction at the fracture surfaces, and (3) a friction model controlling the fracture behavior more specifically. Our approach treats fractures as internal boundary conditions, which naturally allows for different deformation models for the fractured and non-fractured media. The internal boundary conditions that represent fracture surfaces are considered as line pairs (for 2D) or face pairs (for 3D), and the required equations for the fracture behavior are defined between pairs.

To verify the model, the convergence of the method is examined for several benchmark problems. Although the lack of special treatment for the singularity at the fracture tip decreases the convergence rate, we show that the method has approximately 1st-order convergence in displacement for all the cases that we examine.

The approach results in a method with a broad application area. Numerical examples demonstrate the flexibility of the method by considering different models for fracture deformation.

We also show that the method can handle complicated fracture networks and 3D cases. Considering its advantages in the coupling of fluid flow and mechanical deformation, the method is well suited for the further development of simulation tools for a wide range of subsurface engineering applications.

## Acknowledgements

The work was funded by the Research Council of Norway through grant no. 228832/E20 and Statoil ASA through the Akademia agreement. This is a pre-print of an article published in Computational Geosciences. The final authenticated version is available online at: https://doi.org/10.1007/s10596-018-9734-8.

## APPENDIX

The analytical solutions of the induced displacements, $u_x$, $u_y$, and stresses, $\sigma_{xx}$, $\sigma_{yy}$, $\sigma_{xy}$, at any point $(x, y)$ for an infinite two-dimensional homogeneous and isotropic elastic nonporous medium containing a finite small thin fracture with constant normal- and shear-displacement discontinuities are given by Crouch and Starfield [12] as

$$u_x = \Delta u_{n_+}^\Gamma \left(2(1-v)\frac{\partial f}{\partial y} - y\frac{\partial^2 f}{\partial x^2}\right) + \Delta u_{\tau_+}^\Gamma \left(-(1-2v)\frac{\partial f}{\partial x} - y\frac{\partial^2 f}{\partial x \partial y}\right), \quad (39)$$

$$u_y = \Delta u_{n_+}^\Gamma \left((1-2v)\frac{\partial f}{\partial x} - y\frac{\partial^2 f}{\partial x \partial y}\right) + \Delta u_{\tau_+}^\Gamma \left(2(1-v)\frac{\partial f}{\partial y} - y\frac{\partial^2 f}{\partial y^2}\right),$$

and

$$\sigma_{xx} = 2\mu\Delta u_{n_+}^\Gamma \left(2\frac{\partial^2 f}{\partial x \partial y} + y\frac{\partial^3 f}{\partial x \partial y^2}\right) + 2\mu\Delta u_{\tau_+}^\Gamma \left(\frac{\partial^2 f}{\partial y^2} + y\frac{\partial^3 f}{\partial y^3}\right), \quad (40)$$

$$\sigma_{yy} = 2\mu\Delta u_{n_+}^\Gamma \left(-y\frac{\partial^3 f}{\partial x \partial y^2}\right) + 2\mu\Delta u_{\tau_+}^\Gamma \left(\frac{\partial^2 f}{\partial y^2} - y\frac{\partial^3 f}{\partial y^3}\right),$$

$$\sigma_{xy} = 2\mu\Delta u_{n_+}^\Gamma \left(\frac{\partial^2 f}{\partial y^2} + y\frac{\partial^3 f}{\partial y^3}\right) + 2\mu\Delta u_{\tau_+}^\Gamma \left(-y\frac{\partial^3 f}{\partial x \partial y^2}\right),$$

where $\Delta u_{n_+}^\Gamma$ and $\Delta u_{\tau_+}^\Gamma$ are the displacement discontinuities in the normal and shear directions, respectively, $\mu$ is the shear modulus, $v$ is Poisson's ratio, and $f$ is a function of the position $(x, y)$ of the field point relative to the center of the fracture. Denoting the half radius of the fracture as $a$, $f$ is given as

$$f(x,y) = -\frac{1}{4\pi(1-v)} \left( y \left( \tan^{-1}\frac{y}{x-a} - \tan^{-1}\frac{y}{x+a} \right) \right) \quad (41)$$
$$-(x-a)\ln\sqrt{(x-a)^2 + y^2} + (x+a)\ln\sqrt{(x+a)^2 + y^2}.$$